\documentclass[11pt,a4paper]{amsart}

\usepackage[T1]{fontenc}
\usepackage[utf8]{inputenc}
\usepackage[english]{babel}
\usepackage{lmodern}
\usepackage[a4paper,margin=26mm]{geometry}
\usepackage{amsmath,amssymb,amsthm,mathtools}
\usepackage{microtype}
\usepackage{booktabs,tabularx}
\usepackage[hidelinks]{hyperref}

\newcommand{\N}{\mathbb N}
\newcommand{\RCA}{\mathsf{RCA}_0}
\newcommand{\GEN}{1\text{-}\mathrm{GEN}}
\newcommand{\FIP}{\mathrm{FIP}}
\newcommand{\IP}{2\mathrm{IP}}
\newcommand{\sq}{\preceq}
\newcommand{\res}{\upharpoonright}
\newcommand{\cat}{\mathbin{{}^{\frown}}}
\newcommand{\Sel}{\operatorname{Sel}}
\newcommand{\Cof}{\operatorname{Cof}}
\newcommand{\Hit}{\operatorname{Hit}}
\newcommand{\scope}{\operatorname{scope}}
\newcommand{\up}{\mathord{\uparrow}}
\newcommand{\supp}{\mathcal U}
\newcommand{\A}{\mathcal A}

\theoremstyle{plain}
\newtheorem{theorem}{Theorem}[section]
\newtheorem{lemma}[theorem]{Lemma}
\newtheorem{corollary}[theorem]{Corollary}
\theoremstyle{definition}
\newtheorem{definition}[theorem]{Definition}
\theoremstyle{remark}
\newtheorem*{remark}{Remark}

\title[The hierarchy of intersection principles]{
$\mathsf{RCA}_0$ Proves the Collapse of the
Intersection Principles Hierarchy}

\author{Patrizio Cintioli}
\address{Mathematics Division, School of Science and Technology,
University of Camerino, 62032 Camerino, Italy}
\email{patrizio.cintioli@unicam.it}


\begin{document}

\begin{abstract}
This paper shows that the pairwise intersection principle
($2\mathrm{IP}$) implies, over $\mathsf{RCA}_0$, the existence
of a Cohen $1$-generic real relative to every given set.
Together with the known implications, this yields the
equivalence over $\mathsf{RCA}_0$ of the finite intersection
principle (FIP), $2\mathrm{IP}$, and $1$-GEN. In fact, FIP
is equivalent to $n\mathrm{IP}$ for every fixed $n\geq2$.
These results give affirmative answers to Questions~1
and~2 of Cholak, Downey, and Igusa, and to
Question~9.10.26 of Dzhafarov and Mummert.
The proof uses cumulative, nested supports satisfying
a finite Helly property and a decoding argument based
on a least selected label. Neither $\Sigma^0_2$ induction
nor $\Sigma^0_2$ bounding is required.
\end{abstract}

\maketitle

\section{Introduction and conventions}

The \emph{finite intersection principle} (FIP), studied
in reverse mathematics by Dzhafarov and
Mummert~\cite{DM13}, asserts that every nontrivial
countable family of subsets of $\N$ has a maximal
subfamily in which every finite collection of at least
two distinct members has nonempty intersection.

For each fixed $n\geq2$, the \emph{$n$-intersection
principle} ($n\mathrm{IP}$) requires a maximal subfamily
in which every collection of exactly $n$ distinct
members has nonempty intersection. In particular,
$2\mathrm{IP}$ is the pairwise intersection principle.
These are the conventions of the original paper
of Dzhafarov and Mummert~\cite{DM13},
where distinctness is understood extensionally.
The formulation in
\cite[Definition~9.10.17]{DMbook}
instead uses distinct indices and, for FIP,
also requires nonempty singleton intersections.

The family constructed below consists entirely of
nonempty, pairwise distinct sets. Thus, for the
constructed instance, the two formulations of
FIP and $2\mathrm{IP}$ agree.
The known implications involving higher
intersection principles are used with the
conventions of their cited sources.

Earlier degree-theoretic work on FIP and its variants includes
Diamondstone, Downey, Greenberg, and Turetsky~\cite{DDGT}.
Cholak, Downey, and Igusa~\cite[Theorem~5.8]{CDI} proved that
FIP is equivalent, over $\mathsf{RCA}_0$, to the principle
$1$-GEN asserting the existence of a Cohen $1$-generic real
relative to each given set. Their proof that $2\mathrm{IP}$
implies $1$-GEN, however, was formalized over the stronger
system $\mathsf{RCA}_0+I\Sigma^0_2$; see
\cite[Proposition~5.7]{CDI}. They explicitly asked in
\cite[Question~1]{CDI} whether $\mathsf{RCA}_0$ proves
$\mathrm{FIP}\leftrightarrow2\mathrm{IP}$.
Their Question~2 asks whether there exist distinct integers
$m,n\geq2$ such that
$\mathsf{RCA}_0$ proves
$m\mathrm{IP}\leftrightarrow n\mathrm{IP}$,
and raises the same question over
$\mathsf{RCA}_0+B\Sigma^0_2$.
The closely related question whether
$\mathsf{RCA}_0\vdash2\mathrm{IP}\to1\text{-}\mathrm{GEN}$
was later listed by Dzhafarov and Mummert
in~\cite[Question~9.10.26]{DMbook}.

In a different, higher-order setting,
Sanders~\cite[Remark~3.9]{Sanders}
showed that the uniform versions of FIP and of
$n\mathrm{IP}$, for $n\geq2$, are equivalent
to Feferman's search functional $\mu^2$.
These uniform equivalences do not settle the
corresponding nonuniform equivalences over $\RCA$.

The main result of this paper establishes
$2\mathrm{IP}\to1\text{-}\mathrm{GEN}$
already over $\RCA$, thereby resolving these questions.
In particular, for any two distinct fixed integers
$m,n\geq2$, the principles $m\mathrm{IP}$ and
$n\mathrm{IP}$ are equivalent over $\RCA$,
and hence also over $\RCA+B\Sigma^0_2$.

The proof replaces the priority-based labeling
argument by a cumulative system of upward-closed
supports.
The structural labels, stage-truncated halting
tests, and jump-decoding strategy originate in
Cholak, Downey, and Igusa~\cite{CDI}.
The new ingredient is the cumulative, injury-free
organization of the supports.
A projection lemma and a finite Helly
property make pairwise compatibility sufficient
for finite compatibility. When the path decoded
from a maximal solution fails to be generic,
a least \emph{selected} label allows the jump
of the original instance to be recovered.
The crucial minimum concerns a $\Sigma^0_1$
selection predicate, not a $\Pi^0_2$ cofinality
predicate; thus the entire proof remains
within the base theory.

\begin{theorem}\label{thm:main}
Over $\RCA$,
\[
\IP\longrightarrow\GEN.
\]
Consequently,
\[
\RCA\vdash \FIP\longleftrightarrow\IP
\longleftrightarrow\GEN.
\]
\end{theorem}

\subsection{Conventions and representations}

We work in \(\RCA\) with an arbitrary set parameter
\(X\).
All references to computability and computable enumerability
are relative to \(X\), unless otherwise stated.

We write \(2^{<\N}\) for the finite binary strings,
\(\lambda\) for the empty string, and
\(\sigma\sq\tau\) for the (not necessarily proper) prefix relation.
For a string \(\nu\) and a c.e.\ predicate \(E\)
on strings, let
\[
C_\nu=\{\tau:\nu\sq\tau\},
\qquad
\up E=\{\tau:(\exists\rho\sq\tau)\,\rho\in E\}.
\]
The supports will be upward-closed classes of finite strings.
Thus a nonempty intersection has a finite witness.

\begin{quote}
\small
\noindent\textbf{Essential formal convention.}\enspace
The notation \(W_i,K,V_{i,j}\) and the notation for c.e.\ supports
denote predicates equipped with enumerations.
We do not assume that these classes exist as
second-order sets.
Operations on supports abbreviate
\(\Sigma^0_1(X)\) formulas.
The actual sets to which we shall apply \(\IP\) will be
uniformly computable in \(X\) and will be constructed
in Section~\ref{sec:realization}.
\end{quote}

Fix a universal listing \((W_i)_{i\in\N}\)
of the \(X\)-c.e.\ predicates on strings, with finite,
monotone approximations \(W_{i,s}\), normalized so that
\begin{equation}\label{eq:normalization}
W_{i,s}=\varnothing\quad(i>s),
\qquad
\rho\in W_{i,s}\Longrightarrow|\rho|\le s.
\end{equation}
This normalization is obtained by delaying enumerations
and does not change any \(W_i\).

Let \(K_s=K_s^X\) be the finite monotone approximation
to the jump predicate. We write
\[
j\in K\iff(\exists s)\,j\in K_s.
\]
Fix an index \(j_\infty\) of a program that never halts,
uniformly for every oracle.

For a path \(Y\in2^\N\), we say that \(Y\)
\emph{meets} \(W_i\) if one of its prefixes
belongs to \(W_i\).
We say that \(Y\) \emph{avoids} \(W_i\) if
\[
(\exists n)\,
C_{Y\res n}\cap\up W_i=\varnothing.
\]
This meet-or-avoid alternative is equivalent to the usual
genericity condition.
Indeed, if \(Y\) does not meet \(W_i\), avoidance
is equivalent to the existence of a prefix of \(Y\)
with no extension in \(W_i\).

\section{Cumulative supports}

The central construction assigns to each label an upward-closed
support of binary strings. A finite snapshot records requirements
already witnessed along a scope, while a truncated component
retains the information needed for the later jump decoding.

For every string \(\nu\), we record the positive hits
already visible at stage \(|\nu|\):
\begin{equation}\label{eq:snapshot}
M_\nu=
\{i:(\exists\rho\sq\nu)\,\rho\in W_{i,|\nu|}\}.
\end{equation}
By \eqref{eq:normalization},
\(M_\nu\subseteq\{0,\ldots,|\nu|\}\).
Thus \(M_\nu\) is a finite set with a code
computable in \(X\).
Moreover,
\begin{equation}\label{eq:snapshotmonotone}
\nu\sq\mu\Longrightarrow M_\nu\subseteq M_\mu.
\end{equation}

Define the truncated component
\begin{equation}\label{eq:truncation}
V_{i,j}
=
\{\rho:(\exists s)\,
[\rho\in W_{i,s}\land j\notin K_s]\}.
\end{equation}
We have
\begin{align}
j\notin K
&\Longrightarrow V_{i,j}=W_i;
\label{eq:nonhalt}\\
j\in K_h
&\Longrightarrow
(\forall\rho\in V_{i,j})\,|\rho|<h.
\label{eq:haltbound}
\end{align}
For the second implication, every stage witnessing
\(\rho\in V_{i,j}\) precedes \(h\), by monotonicity of \(K_s\).

We use structural labels \(A_\sigma\), for every
string \(\sigma\), and coding labels
\(B_{e,j,\nu}\), for every \(e,j,\nu\).
Their supports are
\[
\supp(A_\sigma)=C_\sigma
\]
and
\begin{equation}\label{eq:support}
\boxed{
\supp(B_{e,j,\nu})=U_{e,j,\nu}
=
C_\nu\cap
\left(
\bigcup_{\substack{i<e\\i\notin M_\nu}}\up W_i
\ \cup\ T_{e,j,\nu}
\right),
}
\end{equation}
where
\[
T_{e,j,\nu}=
\begin{cases}
\up V_{e,j},&e\notin M_\nu,\\
\varnothing,&e\in M_\nu.
\end{cases}
\]

The \emph{scope}, namely the prefix delimiting the support,
is
\[
\scope(A_\sigma)=\sigma,
\qquad
\scope(B_{e,j,\nu})=\nu.
\]
Every support is upward closed under extension and contained
in the cone of its scope.

All components \(i<e\) are untruncated.
Only the current component \(i=e\) depends
on the stage at which \(j\) enters \(K\).
There are no root-label shifts
or injury actions.

\section{Projection, nesting, and the Helly property}

We now establish the structural properties of the supports.
Projection brings them to a common cone, where nesting turns
pairwise compatibility into a common intersection. The argument
must remain valid for internally finite lists in \(\RCA\).

\begin{lemma}[Projection]\label{lem:projection}
If \(\nu\sq\mu\), then
\[
U_{e,j,\nu}\cap C_\mu=C_\mu
\quad\text{or}\quad
U_{e,j,\nu}\cap C_\mu=U_{e,j,\mu}.
\]
There is a flag \(q(e,j,\nu,\mu)\) decidable in \(X\),
using data up to stage \(|\mu|\), that selects
one of the two representations:
if \(q\) holds, the first identity holds;
if \(q\) fails, the second identity holds.

The flag need not decide whether \(U_{e,j,\mu}\)
itself equals all of \(C_\mu\).
\end{lemma}

\begin{proof}
Put \(u=|\mu|\) and define
\begin{align*}
q(e,j,\nu,\mu)\iff{}&
(\exists i<e)\,[i\in M_\mu\setminus M_\nu]\\
&\lor\Bigl[
e\in M_\mu\setminus M_\nu
\ \land\\[-2pt]
&\hspace{12mm}
(\exists s\le u)(\exists\rho\sq\mu)\,
[\rho\in W_{e,s}\land j\notin K_s]
\Bigr].
\end{align*}

If the first disjunct holds, an untruncated component
\(W_i\) present in the old support already contains
a prefix of \(\mu\).
The old support therefore contains all of \(C_\mu\).

If the second disjunct holds, the same conclusion
follows from the presence of a prefix of \(\mu\)
in \(V_{e,j}\).

Now suppose that \(q\) is false.
For every \(i<e\), the exclusions in the two snapshots
coincide, by \eqref{eq:snapshotmonotone}.
Hence the lower components of the two supports,
restricted to \(C_\mu\), are identical.

If \(e\in M_\nu\), the current component is omitted
from both supports.
If \(e\notin M_\mu\), the same \(V_{e,j}\)
appears in both.

It remains to consider the case
\[
e\in M_\mu\setminus M_\nu.
\]
If \(j\notin K_u\), the prefix of \(\mu\)
witnessing \(e\in M_\mu\) belongs to \(V_{e,j}\)
at stage \(u\), contradicting \(\neg q\).
Thus \(j\in K_u\).

Every generator \(\rho\) of \(V_{e,j}\)
now has length less than \(u\),
by \eqref{eq:haltbound}.
If such a generator were compatible with \(\mu\),
it would be a prefix of \(\mu\);
its enumeration stage, which is less than \(u\),
would make \(q\) true.
Consequently,
\[
C_\mu\cap\up V_{e,j}=\varnothing.
\]
The current component of the old support is empty
on this cone, while that of the new support is omitted.
This yields the second identity.
\end{proof}

\begin{lemma}[Nesting at a fixed scope]\label{lem:nesting}
For every \(\mu\):
\begin{enumerate}
\item if \(e<d\), then
\[
U_{e,j,\mu}\subseteq U_{d,k,\mu}
\qquad\text{for all }j,k;
\]
\item for fixed \(e\), the supports \(U_{e,j,\mu}\)
are linearly ordered by inclusion.
\end{enumerate}
\end{lemma}

\begin{proof}
The first inclusion follows from \(V_{e,j}\subseteq W_e\):
every component of the support with index \(e\)
appears in full in the support with index \(d\),
with the same exclusions imposed by \(M_\mu\).

For the second assertion, if \(e\in M_\mu\),
the supports coincide for all \(j\),
since the current component is omitted.

Suppose therefore that \(e\notin M_\mu\).
If neither \(j\) nor \(k\) ever enters \(K\), then
\[
V_{e,j}=V_{e,k}=W_e.
\]
Otherwise consider the first entry of either of them.
If \(j\) enters first, every stage preceding
its entry also precedes the entry of \(k\),
or \(k\) never enters.
Therefore \(V_{e,j}\subseteq V_{e,k}\).
The opposite case is symmetric.

The first entry stage, when it exists, is a
\(\Sigma^0_1\) least number.
\end{proof}

\begin{lemma}[Finite Helly property in \(\RCA\)]\label{lem:helly}
Every internally finite family of labels
with nonempty supports that intersect pairwise
has nonempty total intersection.
\end{lemma}

\begin{proof}
We may remove repetitions.
The empty case, with the usual convention for
an empty intersection, and the case of a single
label are immediate.

For a family of at least two labels,
their scopes are pairwise comparable,
because the supports are contained in the cones of their scopes.
Choose a scope \(\mu\) of maximum length.

Every support has nonempty intersection with \(C_\mu\):
this is witnessed by its intersection with a label
whose scope is \(\mu\);
for that label itself, nonemptiness suffices.

Restrict all supports to \(C_\mu\).
The structural supports become \(C_\mu\).
For the others, compute the flag from
Lemma~\ref{lem:projection}.
True flags again yield \(C_\mu\);
false flags yield supports \(U_{e,j,\mu}\).

If there are no false flags, we are done.
Otherwise choose the least \(e\)
among the finitely many labels with a false flag,
by bounded search.
Among those with this \(e\), consider the finite list
of their indices \(j\).

If one of these indices enters \(K\),
take the first stage at which at least one
of them enters, and choose an index that enters then.
This is just one \(\Sigma^0_1\) least entry stage,
not a choice of a least element in an inclusion ordering.

If none enters, choose any label
from the list: all the corresponding \(V_{e,j}\)
coincide with \(W_e\).
The possible exclusion \(e\in M_\mu\)
is common to all supports in the list.

In either case, the chosen restricted support
is included in those with the same \(e\),
and, by Lemma~\ref{lem:nesting},
in those with a larger index.
It is also included in the full cones.
It is nonempty by the second paragraph;
a witness in it therefore belongs to every support.

The bounded choices concern decidable finite data.
The only possible unbounded choice is the first entry
of an element of the list into \(K\),
which is permitted by \(I\Sigma^0_1\).
Thus the argument also works for families
of internally nonstandard finite cardinality.
\end{proof}

\section{Realization as a computable family}
\label{sec:realization}

The preceding combinatorial assertions concern enumerable
classes of strings. We now code them as an actual family of
sets: common odd tokens record intersection witnesses, and
private even markers identify the family members. The coding
uses no arithmetical comprehension and excludes empty sets.

Fix codes for all potential labels
and a uniform enumeration of positive assertions
\(\tau\in\supp(L)\).
We use finite monotone approximations \(\supp_s(L)\),
decidable in \(X\), with only finitely many labels
and strings examined at each stage.
Every true membership assertion appears at some stage.

This enumeration is obtained by dovetailing
formulas \eqref{eq:truncation} and \eqref{eq:support},
with restrictions on the codes at the current stage.
The truncation in \(V_{i,j}\) concerns the stage
witnessing \(\rho\in W_{i,s}\), as in the
definition, and not a later stage
at which the witness is processed.

Enumerate without repetitions only those labels
for which a witness to nonemptiness has appeared.
Let
\[
L_0,L_1,L_2,\ldots
\]
be the resulting enumeration.

All structural labels \(A_\sigma\) have nonempty support
and are inserted according to a fixed computable scan.
Consequently, \(\RCA\) proves that each \(L_n\) exists:
the first \(n+1\) distinct structural witnesses give
a finite bound for the appearance of at least
\(n+1\) labels.
The function \(n\mapsto L_n\) is total
and computable in \(X\).

Define a family
\[
\A=(D_n)_{n\in\N}
\]
of subsets of \(\N\).
The set \(D_n\) contains the single even number \(2n\).
For odd numbers, let
\begin{equation}\label{eq:tokens}
2\langle\tau,s\rangle+1\in D_n
\iff
\begin{cases}
L_n\text{ has already been indexed by stage }s,\\
\tau\in\supp_s(L_n).
\end{cases}
\end{equation}
The two conditions on the right are required jointly.

This test is decidable in \(X\):
one need only simulate the finite stage \(s\).
If fewer than \(n+1\) labels have appeared by that stage,
the answer is negative.
Thus the relation \(z\in D_n\) defines
a set by \(\Delta^0_1(X)\) comprehension.

All the \(D_n\) are nonempty and pairwise distinct,
thanks to the private markers \(2n\).
For two distinct indices, a common element
must be odd; hence
\begin{equation}\label{eq:intersectioncoding}
D_n\cap D_m\ne\varnothing
\iff
\supp(L_n)\cap\supp(L_m)\ne\varnothing.
\end{equation}

The same equivalence holds for every finite family
of at least two distinct indices.
In one direction, a common odd element gives
a witness \(\tau\).
Conversely, given a common witness \(\tau\)
for the supports, \(B\Sigma^0_1\), available in \(\RCA\),
collects the finitely many required enumeration
and indexing stages.
A later common stage provides the shared odd number
in \eqref{eq:tokens}.

For a single label, both sides of the equivalence
are true by construction.
Thus the Helly property transfers
to the actual family \(\A\).

Applying \(\IP\) to \(\A\), we obtain a maximal subfamily.
We may represent it by a total function \(f\)
such that
\[
(D_{f(r)})_{r\in\N}
\]
is maximal with the pairwise intersection property.

If the formulation of the principle instead returns
a sequence of sets \((E_r)\),
we recover \(f(r)\) by searching for the unique \(n\)
with \(2n\in E_r\).
The search terminates because each \(E_r\) is a member
of \(\A\).
We do not decide equality between arbitrary sets.

For a potential label \(L\), define
\begin{equation}\label{eq:selected}
\Sel(L)\iff(\exists r)\,L_{f(r)}=L.
\end{equation}
This is a \(\Sigma^0_1(X\oplus f)\) property.
We do not assume that the range of \(f\)
exists as a set.
A label with empty support is not indexed
and is not selected.

By \eqref{eq:intersectioncoding}, two selected labels
have intersecting supports.
Conversely, if a label has nonempty support
and is compatible with every selected label,
its set can be added to the subfamily
while preserving \(2\)IP.
Maximality therefore forces it to be selected.

Formally, to add \(D_n\) one merely prepends it
to the sequence \((D_{f(r)})_r\).
This operation is available in \(\RCA\)
and uses the extensional maximality of the subfamily.

\section{The path and the characterization of selection}

The structural labels selected by a maximal pairwise-intersecting
subfamily determine an infinite path through the binary tree.
We then characterize the selection of any label by the cofinality
of its support along that path.

\begin{lemma}[Structural path]\label{lem:path}
From \(X\oplus f\) one can compute a unique path
\(Y\in2^\N\) such that
\[
\Sel(A_\sigma)\iff\sigma\sq Y.
\]
The argument can be formalized in \(\RCA\).
\end{lemma}

\begin{proof}
The support of \(A_\lambda\) intersects every nonempty support.
Therefore \(A_\lambda\) is selected by maximality.
The scopes of selected structural labels are
pairwise comparable.

Suppose that \(A_\sigma\) is selected and neither
of its two children \(A_{\sigma\cat0}\),
\(A_{\sigma\cat1}\) is selected.
By maximality, there are selected labels
\(H_0,H_1\) with
\[
\supp(H_b)\cap C_{\sigma\cat b}=\varnothing
\qquad(b=0,1).
\]
Indeed, if a child had no obstruction
among the selected labels,
it could be added to the solution.

The three labels \(A_\sigma,H_0,H_1\)
have nonempty total intersection,
by Lemma~\ref{lem:helly}.
Let \(\tau\) be a common witness.

Since the supports are upward closed,
we can extend \(\tau\) to a string
of length greater than \(|\sigma|\).
That extension belongs to the cone of one
of the two children, contradicting the corresponding
disjointness condition.

Thus at least one child is selected.
They cannot both be selected,
since their cones are disjoint.

The formula
\[
(\exists\sigma)\,
[|\sigma|=n\land\Sel(A_\sigma)]
\]
is \(\Sigma^0_1(X\oplus f)\).
\(\Sigma^0_1\) induction gives a selected node
at every level.
It is unique, and the nodes are compatible.

By searching the enumeration for the node
at the required level, we compute the bits of \(Y\).
The totality just proved permits
\(\Delta^0_1(X\oplus f)\) comprehension.
Uniqueness and comparability give both directions
of the equivalence in the statement.
\end{proof}

\begin{definition}
A support \(\supp(L)\) is \emph{cofinal along \(Y\)} if
\begin{equation}\label{eq:cofinal}
\Cof_Y(L)
\iff
(\forall n)(\exists\tau\succeq Y\res n)\,
\tau\in\supp(L).
\end{equation}
This is a \(\Pi^0_2(X\oplus Y)\) property.
\end{definition}

\begin{lemma}[Cofinality and maximality]
\label{lem:cofinalselected}
For every potential label \(L\),
\begin{equation}\label{eq:cofinalselected}
\boxed{\Sel(L)\iff\Cof_Y(L).}
\end{equation}
\end{lemma}

\begin{proof}
If \(L\) is selected, its support intersects
every \(C_{Y\res n}\), since all the
\(A_{Y\res n}\) are selected.

Conversely, suppose \(\Cof_Y(L)\).
Its scope is a prefix of \(Y\):
if it were incompatible with a prefix of \(Y\),
it could not be cofinal.

The same observation applies to every selected label
\(H\), by the direction already proved.
Let \(\mu\) be the longer of the two scopes.
The restrictions of \(\supp(L)\) and \(\supp(H)\)
to \(C_\mu\) are both nonempty.

By projection and nesting, these restrictions
are comparable by inclusion,
including the case of the full cone.
Thus they intersect.

The label \(L\) has nonempty support
and is therefore indexed.
It is compatible with every selected label,
and maximality forces it into the solution.

This also deals with potential labels
that are not indexed:
for them both properties are false.
\end{proof}

\begin{remark}
Two arbitrary open sets may both be cofinal
along the same path without intersecting.
Here the implication is justified by projection
to a common scope and nesting.
This is a specific property of the family
\eqref{eq:support}.
\end{remark}

\section{Stabilization of positive hits only}

\begin{lemma}[Exact finite snapshot]
\label{lem:stabilization}
For every \(a\) and every path \(Y\),
there exists \(t>a\) such that, putting \(\nu=Y\res t\),
we have, for each \(i\le a\),
\begin{equation}\label{eq:stablesnapshot}
i\in M_\nu\iff\Hit(i,Y),
\qquad
\Hit(i,Y)\equiv
(\exists\rho\sq Y)\,\rho\in W_i.
\end{equation}
This is provable in \(\RCA\).
\end{lemma}

\begin{proof}
Let
\[
c(s)=
\#\{i\le a:
(\exists\rho\sq Y\res s)\,\rho\in W_{i,s}\}.
\]
The function \(c\) is computable in \(X\oplus Y\),
nondecreasing and bounded by \(a+1\).

The predicate
\[
P(r)\equiv(\exists s)\,[a+1-c(s)=r]
\]
is \(\Sigma^0_1(X\oplus Y)\) and nonempty.
The \(\Sigma^0_1\) least number principle,
available from \(I\Sigma^0_1\),
gives its least element \(r_0\).

A stage \(s_0\) witnessing \(P(r_0)\)
achieves the maximum of \(c\).
By monotonicity,
\[
c(s)=c(s_0)\qquad(s\ge s_0).
\]
No new index \(i\le a\) can therefore
be recorded subsequently.

Choose \(t>\max\{a,s_0\}\).
If \(W_i\) contains a prefix \(\rho\) of \(Y\),
then at a sufficiently late stage that prefix
appears in \(W_{i,s}\) and is a prefix of \(Y\res s\).
Thus index \(i\) must already have been present
in the snapshot at \(t\).

The reverse implication in \eqref{eq:stablesnapshot}
is immediate.
\end{proof}

\begin{remark}
A bounded sequence of positive events is stabilized:
each index \(i\le a\) can contribute just one new event.
We do not collect a family of bounds for components
that are individually noncofinal.
The latter operation will be replaced by a single
bound for a cumulative support.
\end{remark}

\section{The least selected index and decoding the jump}

Suppose that the path \(Y\) from
Lemma~\ref{lem:path} is not \(1\)-generic
relative to \(X\).
Then there is an index \(a\) such that \(Y\)
neither meets nor avoids \(W_a\).
In particular,
\begin{equation}\label{eq:badrequirement}
\neg\Hit(a,Y),
\qquad
(\forall n)\,
C_{Y\res n}\cap\up W_a\ne\varnothing.
\end{equation}

Apply Lemma~\ref{lem:stabilization}
through index \(a\) and fix \(\nu=Y\res t\)
with an exact snapshot.

Every residual component \(W_i\) with
\(i\le a\) and \(i\notin M_\nu\)
therefore misses every prefix of \(Y\).

For \(r\le a\), define the potential labels
\begin{equation}\label{eq:N}
N_r=B_{r,j_\infty,\nu}.
\end{equation}
Since \(j_\infty\notin K\), their supports are
\begin{equation}\label{eq:Nsupport}
\supp(N_r)
=
C_\nu\cap
\bigcup_{\substack{i\le r\\i\notin M_\nu}}\up W_i.
\end{equation}

By \eqref{eq:badrequirement},
\(a\notin M_\nu\) and \(\supp(N_a)\)
is cofinal along \(Y\).
By Lemma~\ref{lem:cofinalselected},
\(N_a\) is selected.

Therefore choose
\begin{equation}\label{eq:minselected}
\boxed{
e=\min\{r\le a:\Sel(N_r)\}.
}
\end{equation}
This choice uses the least number principle for a
\(\Sigma^0_1(X\oplus f)\) predicate, by \eqref{eq:selected}.
We are not taking the least index satisfying
the \(\Pi^0_2\) cofinality property of the individual \(W_r\).

Now let
\begin{equation}\label{eq:baseline}
\mathcal R
=
C_\nu\cap
\bigcup_{\substack{i<e\\i\notin M_\nu}}\up W_i.
\end{equation}
If \(e=0\), then \(\mathcal R=\varnothing\).
If \(e>0\), we have the exact identity
\begin{equation}\label{eq:baselineidentity}
\mathcal R=\supp(N_{e-1}).
\end{equation}

In the latter case, where \(e>0\),
minimality gives \(\neg\Sel(N_{e-1})\).
Lemma~\ref{lem:cofinalselected}
then implies that \(\mathcal R\) is not cofinal.

For clarity, note also that
\begin{equation}\label{eq:eoutsidesnapshot}
e\notin M_\nu.
\end{equation}
Indeed, if $e=0$ and $e\in M_\nu$, then
$\supp(N_e)=\varnothing$, contradicting its selection.
If $e>0$ and $e\in M_\nu$, the current component of
$N_e$ is omitted, and hence
$\supp(N_e)=\supp(N_{e-1})$.
The equivalence between cofinality and selection would
then contradict the minimality of $e$.

In either case, there exists a \emph{single} level \(b\)
such that
\begin{equation}\label{eq:onebound}
C_{Y\res b}\cap\mathcal R=\varnothing.
\end{equation}
This is an existence assertion, not a claim of uniqueness
of $b$, and no bound is computed from the solution.
In particular, no separate bound is chosen for each
\(i<e\).

\begin{lemma}[A column decoding \(K\)]
\label{lem:decoder}
For the \(e,\nu\) just chosen and every \(j\),
\begin{equation}\label{eq:decoder}
\boxed{
\Sel(B_{e,j,\nu})\iff j\notin K.
}
\end{equation}
\end{lemma}

\begin{proof}
If \(j\notin K\), then by \eqref{eq:nonhalt}
and \eqref{eq:Nsupport},
\[
U_{e,j,\nu}=\supp(N_e).
\]
This support is cofinal because \(N_e\)
is selected.
Therefore \(B_{e,j,\nu}\) is selected
by Lemma~\ref{lem:cofinalselected}.

Now suppose \(j\in K_h\).
The support \(U_{e,j,\nu}\) consists of the part
\(\mathcal R\) and the current component,
when present.

The latter is generated by strings
\(\rho\in V_{e,j}\subseteq W_e\)
of length less than \(h\).
If this component is present, then \(e\notin M_\nu\),
and the exact snapshot guarantees that none
of these strings is a prefix of \(Y\).

Each \(\rho\) is therefore incompatible
with \(Y\res(h+1)\):
since its length is less than \(h\),
if it were compatible it would be a prefix.
Consequently,
\[
C_{Y\res(h+1)}\cap T_{e,j,\nu}=\varnothing.
\]
If the component is omitted,
the same equality is obvious.

Together with \eqref{eq:onebound}, this gives
\[
C_{Y\res\max\{b,h+1,t\}}\cap U_{e,j,\nu}
=
\varnothing.
\]
The label is not cofinal
and therefore is not selected.
\end{proof}

By \eqref{eq:decoder}, the jump predicate
has two complementary definitions:
\[
j\in K\iff(\exists s)\,j\in K_s,
\qquad
j\notin K\iff\Sel(B_{e,j,\nu}).
\]
Both right-hand sides are
\(\Sigma^0_1(X\oplus f)\), with the numerical
parameters \(e,\nu\) fixed.

Thus \(\Delta^0_1\) comprehension produces
a set \(J\) with
\begin{equation}\label{eq:jumpexists}
j\in J\iff(\exists s)\,j\in K_s.
\end{equation}
This is exactly \(X'\).

Computationally, one alternates between
searching for the entry of \(j\) into \(K\)
and searching for label \(B_{e,j,\nu}\)
in the enumeration of the solution.
Exactly one of these two searches terminates.

The decoder need not know or compute
the bound \(b\): that bound is used only
in the correctness proof.

\section{Proof of the main theorem and consequences}

\begin{lemma}\label{lem:jumpgeneric}
In \(\RCA\), if the jump \(X'\)
exists as a set, then there exists a \(1\)-generic
relative to \(X\).
\end{lemma}

\begin{proof}
Construct an increasing sequence of strings
\((\sigma_s)_s\), with \(\sigma_0=\lambda\).
At stage \(s\), use \(X'\) to decide
the \(\Sigma^0_1(X)\) question
\[
(\exists\tau\succeq\sigma_s)
(\exists\rho\sq\tau)\,\rho\in W_s.
\]

If the answer is positive, search for witnesses
\(\tau,\rho\) and put
\[
\sigma_{s+1}=\tau\cat0.
\]
If the answer is negative, put
\[
\sigma_{s+1}=\sigma_s\cat0.
\]

The jump decides the question uniformly,
using the index of an \(X\)-oracle program
that searches for witnesses.
The search on the positive branch terminates
by correctness of the decision.

Totality of the construction at all stages
is proved by \(\Sigma^0_1(X\oplus X')\) induction
on the existence of finite computation traces.
There is no need to include in the inductive formula
the assertion that all preceding requirements
have been satisfied.

Since \(|\sigma_s|\ge s\),
the union of the strings defines a set \(G\)
by \(\Delta^0_1(X\oplus X')\) comprehension.

At stage \(s\), the positive branch meets \(W_s\),
since \(\rho\sq G\).
The negative branch certifies
\[
C_{\sigma_s}\cap\up W_s=\varnothing.
\]
Therefore \(G\) meets or avoids every \(W_s\)
and is thus \(1\)-generic relative to \(X\).
\end{proof}

\begin{proof}[Proof of Theorem~\ref{thm:main}]
Let \(X\) be arbitrary.
Construct the family computable in \(X\)
from Section~\ref{sec:realization}
and apply \(\IP\).

From the maximal solution we obtain \(f\),
and from Lemma~\ref{lem:path} we obtain the path \(Y\).

If \(Y\) is \(1\)-generic relative to \(X\),
we already have the desired witness.

If it is not, choose a requirement \(W_a\)
that \(Y\) neither meets nor avoids.
Lemma~\ref{lem:stabilization}
and the construction in Section 7
produce the set \(X'\) in \(\RCA\).
Lemma~\ref{lem:jumpgeneric}
then produces a \(1\)-generic relative to \(X\).

In either case, a witness to \(\GEN\) exists.
Since \(X\) was arbitrary,
\[
\RCA\vdash\IP\to\GEN.
\]

Proposition~5.4 of~\cite{CDI} gives
$\RCA\vdash\GEN\to\FIP$, while their
Proposition~5.6 gives $\RCA\vdash\FIP\to\IP$.
Together with the implication just proved, these form
a cycle, and therefore
\[
\RCA\vdash\FIP\longleftrightarrow\IP
\longleftrightarrow\GEN.
\]
\end{proof}

\begin{corollary}\label{cor:alln}
For each fixed integer $n\geq2$,
$\RCA$ proves that $n\mathrm{IP}$, $\IP$, $\FIP$,
and $\GEN$ are equivalent. In particular, for every
fixed $m,n\geq2$,
\[
\RCA\vdash m\mathrm{IP}\longleftrightarrow n\mathrm{IP}.
\]
\end{corollary}

\begin{proof}
Combine Theorem~\ref{thm:main}
with the known implications
\[
\FIP\to n\mathrm{IP}\to\IP,
\qquad
\FIP\leftrightarrow\GEN,
\]
recalled in
\cite[Propositions 5.4--5.6]{CDI}.
\end{proof}

\section{Formalization in the base theory}

The proof is intended literally as a derivation in
$\RCA$, not merely as an argument valid in its
$\omega$-models. For reference, the following table
identifies the resources involved at the potentially
sensitive steps.

\begin{center}
\small
\renewcommand{\arraystretch}{1.2}
\begin{tabularx}{\textwidth}{@{}p{.39\textwidth}X@{}}
\toprule
Step & Resource in \(\RCA\)\\
\midrule
Snapshot \(M_\nu\) and projection flag
&
Finite data decidable in \(X\).
\\

Actual family \((D_n)\)
&
\(\Delta^0_1(X)\) comprehension and stage coding.
\\

Helly property for an internally finite list
&
Bounded search and, if needed,
a least \(\Sigma^0_1\) entry stage.
\\

Witnesses for shared support elements
&
\(B\Sigma^0_1\) bounding.
\\

Path from \(f\)
&
\(\Sigma^0_1(X\oplus f)\) induction
and \(\Delta^0_1\) comprehension.
\\

Stabilization up to \(a\)
&
\(\Sigma^0_1\) minimization of the reversed values
of a bounded counting function.
\\

Choice of \(e\)
&
\(\Sigma^0_1(X\oplus f)\) minimization
for selection.
\\

Lower bound \(b\)
&
Noncofinality of one support,
\(\supp(N_{e-1})\), when \(e>0\).
\\

Existence of \(X'\)
&
Two complementary enumerations
and \(\Delta^0_1\) comprehension.
\\

Generic from \(X'\)
&
Computable recursion and
\(\Sigma^0_1(X\oplus X')\) induction.
\\
\bottomrule
\end{tabularx}
\end{center}

The minimum in \eqref{eq:minselected}
uses information provided by the solution:
the positive appearance of a label
in its enumeration.

The preceding support \(N_{e-1}\)
collects all lower components into one class.
Its nonselection gives a single avoidance level.
This is the mechanism that replaces
the choice of a \(\Pi^0_2\) minimum
and the collection of bounds for lower components.

The choices of \(a,t,e,\nu\) are existential
and nonuniform.
We do not assert the existence of a single algorithm
that determines, from a solution, which of the two cases
\emph{\(Y\) is generic} or \emph{\(X'\) is decodable}
holds, or that uniformly finds the parameters
of the decoder.

Thus the conclusion concerns the implication
over \(\RCA\).
It does not automatically provide a Weihrauch reduction
from \(\GEN\) to \(\FIP\).

\section{Concluding remarks}

Theorem~\ref{thm:main} answers Question~1 of
Cholak--Downey--Igusa~\cite{CDI} affirmatively: the
pairwise and finite intersection principles are equivalent
already over $\RCA$. Corollary~\ref{cor:alln} also answers
their Question~2, in its $\RCA$ formulation, for every
pair of fixed integers $m,n\geq2$.
Since these equivalences are provable already
in $\RCA$, they also hold over
$\RCA+B\Sigma^0_2$, answering the second part
of Question~2 as well.
Moreover, the
$\Sigma^0_2$-induction hypothesis in the known derivation
of $2\mathrm{IP}\to 1$-GEN is not needed for the implication.

The argument relies on two features of the new supports:
their finite Helly property, which is available for
internally finite lists, and the cumulative lower part
of a least selected label, which admits a single avoidance
level. The actual family is $\Delta^0_1(X)$, and its
private markers allow solutions to be decoded without
requiring the range of an enumerating function as a
second-order set.

The proof is nonuniform: the choice between an
already generic path and a decodable jump, as well
as the parameters of the jump decoder, is
existential. The equivalence over $\RCA$ does not
by itself establish a Weihrauch reduction between
the corresponding problems.

In fact, a companion manuscript~\cite{CintioliW}
establishes
\[
1\text{-}\mathrm{GEN}\not\leq_{\mathrm W}\mathrm{FIP}.
\]
Thus, although $\RCA$ proves that FIP implies
$1$-GEN, this implication cannot be realized
by a Weihrauch reduction in the representations
considered there.

\section*{Acknowledgments}

The main result was discovered by ChatGPT 6 Pro (OpenAI),
which was also used in developing the proof and preparing
the manuscript. A preliminary version of the proof was
reviewed using Claude Opus 5.5 (Anthropic), and a later
version underwent a separate mathematical audit using
Astra 6.1 Ultra (OpenAI).
The author takes full responsibility for the correctness
of the mathematical statements and proofs, for the
accuracy of the citations and attributions, and for
the final content of the manuscript.

\end{document}